\documentclass[12pt,a4paper,twoside]{amsart}
\usepackage{a4}

\usepackage[cp1251]{inputenc}
\usepackage[T2A]{fontenc}

\usepackage[russian, english]{babel}
\usepackage{amsmath, amsfonts, amssymb, amsthm}

\usepackage{pgf,tikz}
\usetikzlibrary{arrows}

\newtheorem{theorem}{Theorem}

\newtheorem*{bla}{Blaschke's Rolling Theorem}

\theoremstyle{remark}
\newtheorem{remark}{Remark}

\theoremstyle{definition}

\newtheorem*{ack}{Acknowledgments}

\DeclareMathOperator{\normc}{n}

\DeclareMathOperator{\dist}{dist}
\DeclareMathOperator{\grad}{grad}

\DeclareMathOperator{\sn}{sn}

\numberwithin{equation}{section}

\title[Comparison theorem for support functions]{Comparison theorem for support functions of hypersurfaces}
\author[A.~Borisenko]{Alexander Borisenko}

\address[A.~Borisenko]{Department of Mathematical Analysis and Optimization, Sumy State University\\
 Rimskogo~- Korsakova str. 2, 40007, Sumy, Ukraine}

\email{aborisenk@gmail.com}

\author[K.~Drach]{Kostiantyn Drach}

\address[K.~Drach]{Geometry Department \\V.N. Karazin Kharkiv National University\\Svobody Sq. 4, 61022, Kharkiv\\Ukraine}

\address{Department of Mathematical Analysis and Optimization \\ Sumy State University\\ Rimskogo~- Korsakova str. 2, 40007, Sumy \\ Ukraine}
		
\email{drach@karazin.ua, kostya.drach@gmail.com}       



\begin{document}

\maketitle

\begin{abstract}
For a convex domain $D$ that is enclosed by the hypersurface $\partial D$ of bounded normal curvature, we prove an angle comparison theorem for angles between $\partial D$ and geodesic rays starting from some fixed point in $D$, and the corresponding angles for hypersurfaces of constant normal curvature. Also, we obtain a comparison theorem for support functions of such surfaces. As a corollary, we present a proof of Blaschke's Rolling Theorem.
\end{abstract}

\section{Preliminaries and the main results}

Is it known the following theorem due to W.~Blaschke: 
\begin{bla}
\label{blth}
Let $\mathbb M^{m}(c)$ be an $m$-dimensional space of constant curvature equal to $c$, $D \subset \mathbb M^{m}(c)$ be a convex body with the $C^r$-smooth boundary $\partial D$ ($r\geqslant 2$), and $P \in \partial D$ be an arbitrary point. Let $\partial D_\lambda \subset \mathbb M^m(c)$ be a complete hypersurface of constant normal curvature equal to some $\lambda > 0$, and suppose that $\partial D_\lambda $ touches $\partial D$ at $P$ so that their inner unit normals coincide.

A. If normal curvatures $k_{\normc}$ of the hypersurface $\partial D$ at all points and in all directions satisfy the inequality $k_{\normc} \geqslant \lambda$, then $\partial D$ lies entirely in the closed convex domain bounded by $\partial D_\lambda$. 

B. If normal curvatures of the hypersurface $\partial D$ at all points and in all directions satisfy the inequality $\lambda \geqslant k_{\normc}$, then the hypersurface $\partial D_\lambda$ lies in $D$. 

Moreover, the hypersurfaces $\partial D$ and $\partial D_\lambda$ can intersect only by a domain that contains the point $P$.
\end{bla}

For the Euclidean space this theorem was first proved in~\cite{Bla}; for the general case of constant curvature spaces see~\cite{Kar, M, How}.

It appears that Blaschke's Rolling Theorem can be obtained as a corollary from the following comparison theorems for angles between the radius-vector of a hypersurface and its normals. In order to give exact statements, we need to agree on some notations.

Everywhere below let $\mathbb M^{m}$ be a complete simply-connected $m$-dimensional Riemannian manifold such that its sectional curvatures $K_\sigma$ in a direction of a $2$-plane $\sigma \subset T\mathbb M^m$ satisfy the inequality $c_2 \geqslant K_\sigma \geqslant c_1$ with some constants $c_1$ and $c_2$. Furthermore, let $D \subset \mathbb M^{m}$ be a closed domain with the boundary $\partial D$ being a $C^r$-smooth hypersurface ($r\geqslant 2$). For $c_2 > 0$ we will additionally assume that the domain $D$ lies inside a geodesic sphere of radius $\pi/(2\sqrt{c_2})$. 

By $t_Q(\cdot)=\dist(Q,\cdot)$ denote a distance function from some point $Q \in D$ defined on $\mathbb M^{m}\backslash \{Q\}$, and let $\partial_{t_Q}$ be a gradient vector field of the function $t_Q$, and $\rho_Q$ be a restriction of $t_Q$ on $\partial D$: $\rho_Q(\cdot) = t_Q(\cdot)\left|_{\partial D}\right.$.

\begin{theorem}
\label{anglecompth1}

Suppose $D \subset \mathbb M^m$ and $D_{k_1} \subset \mathbb M^m (c_1)$ are closed domains such that normal curvatures $k_{\normc}$ of the hypersurface $\partial D$ at any point and in any direction with respect to the inner unit normal field $N$ satisfy the inequality $$k_{\normc} \geqslant k_1 > 0,$$ and normal curvatures of $\partial D_{k_1}$ are constant and equal to $k_1$ with respect to the inner unit normal field $N_1$. 
Let $O \in D$ and $O_1 \in D_{k_1}$ be points with $\dist(O,\partial D) = \dist(O_1, \partial D_{k_1})$; then at all points $P \in \partial D$ and $P_1 \in \partial D_{k_1}$ such that $$\rho_{O}(P) = \rho_{O_1}(P_1),$$ the inequality
\begin{equation}
\label{actheq1}
\left|\left<N,\partial_{t_O}\right>\right| (P) \geqslant \left|\left<N_1,\partial_{t_{O_1}}\right>\right|(P_1)
\end{equation}
holds.
\end{theorem}

Recall that a function $h_Q \colon \partial D \to (0,+\infty)$ defined as $$h_Q = \rho_Q \cdot |\left<N, \partial_{t_Q}\right>|$$ is called a \textit{support function} of the hypersurface $\partial D \subset \mathbb M^m$ with respect to a point $Q \in D$ (see~\cite[chapter 6, \S 5]{BZ}). 

Using Theorem~\ref{anglecompth1} we can obtain the following comparison theorem for support functions.

\begin{theorem}
\label{supfunccompth1}
Let $D \subset \mathbb M^m$ and $D_{k_1} \subset \mathbb M^m (c_1)$ be closed domains such that normal curvatures $k_{\normc}$ of the hypersurface $\partial D$ satisfy the inequality $$k_{\normc} \geqslant k_1 > 0,$$ and normal curvatures of $\partial D_{k_1}$ are constant and equal to $k_1$.
Let $O \in D$ and $O_1 \in D_{k_1}$ be points with $\dist(O,\partial D) = \dist(O_1, \partial D_{k_1})$; then at all points $P \in \partial D$ and $P_1 \in \partial D_{k_1}$ such that $\rho_{O}(P) = \rho_{O_1}(P_1)$, the inequality
\begin{equation*}
\label{sfctheq1}
h_O(P) \geqslant h_{O_1} (P_1).
\end{equation*}
holds.
\end{theorem}

For Theorems~\ref{anglecompth1} and~\ref{supfunccompth1} also holds the following dual result.

\begin{theorem}
\label{asfcompth}
Suppose $D \subset \mathbb M^m$ and $D_{k_2} \subset \mathbb M^m (c_2)$ are closed domains such that normal curvatures $k_{\normc}$ of the hypersurface $\partial D$ with respect to the inner unit normal field $N$ satisfy the inequality $$k_2 \geqslant k_{\normc} > 0,$$ and normal curvatures of $\partial D_{k_2}$ are constant and equal to $k_2$ with respect to the inner unit normal field $N_2$. Let $O \in D$ and $O_2 \in D_{k_2}$ be points with $\dist(O,\partial D) = \dist(O_2, \partial D_{k_2})$; then at all points $P \in \partial D$ and $P_2 \in \partial D_{k_2}$ for which the distances $\rho_{O}(P)$ and $\rho_{O_2}(P_2)$ are equal, the inequalities
\begin{equation*}
\label{actheq2}
\begin{aligned}
\left|\left<N_2,\partial_{t_{O_2}}\right>\right| &\geqslant \left|\left<N,\partial_t\right>\right|,\\
h_{O_2} &\geqslant h_{O}
\end{aligned}
\end{equation*}
hold.
\end{theorem}

\begin{remark}
Actually, in Theorem~\ref{asfcompth} we need only the weaker restriction $c_2 \geqslant K_\sigma$ on sectional curvatures of the manifold $\mathbb M^{m}$.
\end{remark}

\begin{remark}
Theorems~\ref{anglecompth1}~-- \ref{asfcompth} will remain true if we replace the convex domain $D$ with a star-shaped domain of normal curvatures bounded above or below by a non-zero number $\lambda$.
\end{remark}

\section{Proof of Theorem~\ref{anglecompth1}}

In this section we will prove Theorem~\ref{anglecompth1} using the similar technique as in~\cite{BorDr2}, but our proof will be shorter.

Let $Q \in \partial D$ and $Q_1 \in \partial D_{k_1}$ be points such that $\dist(O, \partial D) = t_O(Q)$ and $\dist(O_1, \partial D_{k_1}) = t_{O_1} (Q_1)$. By $d$ denote the distance $t_O(Q) =  t_{O_1} (Q_1)$. Observe that inequality~(\ref{actheq1}) holds at $Q$ and $Q_1$.

In the manifolds $\mathbb M^{m}$ and $\mathbb M^{m}(c_1)$ let us introduce polar coordinate systems with origins, respectively, at $O$ and $O_1$. By hypothesis of the theorem, both hypersurfaces lie in the regularity regions of these systems of coordinates. Moreover, since the second fundamental forms of $\partial D$ and $\partial D_{k_1}$ are positively defined, the hypersurfaces bound the convex regions. Thus they both can be explicitly defined in the introduced coordinate systems. 

Suppose $\gamma(t)$ and $\gamma_1(t)$ are integral trajectories of the gradient vector fields for the functions $\rho_O$ and $\rho_{O_1}$ passing through the points $P$ and $P_1$, and parametrized by a parameter $t$ measuring the distance from the corresponding origin. We note that $Q$ and $Q_1$ are limit points of, respectively, $\gamma$ and $\gamma_1$, and $\gamma(d) = Q$, $\gamma_1(d)=Q_1$. It appears that along these integral trajectories the following equalities hold (see ~\cite{Bor1} for details)
\begin{equation}
\label{normcF}
k_{\normc}(t) = \left|\left<N,\partial_{t_O}\right> \right| (t) \cdot \mu_{\normc}(t) + \frac{d}{dt}\left|\left<N,\partial_{t_{O}}\right>\right|,
\end{equation}
\begin{equation}
\label{normcF1}
k_1 = \left|\left<N_1,\partial_{t_{O_1}}\right> \right| (t) \cdot \mu_{\normc}^{c_1} (t) + \frac{d}{dt}\left|\left<N_1,\partial_{t_{O_1}}\right>\right|,
\end{equation}
where $\mu_{\normc}^{c_1}(t)$ is the normal curvature of a sphere of radius $t$ in $\mathbb M^{m}(c_1)$; $k_{\normc}(t)$ is the normal curvature of $\partial D$ taken at the point $\gamma(t)$ in the direction of the vector $\dot \gamma(t)$; $\mu_{\normc}(t)$ is the normal curvature of the geodesic sphere $S^{m-1} \subset \mathbb M^m$ of radius $t$ and center $O$ taken at the point $\gamma(t)$ in the directions of the projection of $\dot \gamma(t)$ on the tangent space $T_{\gamma(t)}S^{m-1}$. All normal curvatures are calculated with respect to the corresponding inner normal vector fields.

It is known that $\mu_{\normc}^{c_1}(t) = {\sn_{c_1}'(t)}/{\sn_{c_1}(t)},$ where
\begin{equation*}
\label{sin}
\sn_{c_1}(t) = \left\{
\begin{aligned}
&\frac{1}{\sqrt{c_1}}\sin \sqrt{c_1} t, \text{ for }c_1 > 0 \\
& t, \text{ for }c_1 = 0 \\
&\frac{1}{\sqrt{-c_1}}\sinh \sqrt{-c_1} t, \text{ for }c_1 < 0. \\
\end{aligned}
\right.
\end{equation*}

By the comparison theorem for normal curvatures of spheres (see~\cite[chapter 6, \S 5]{Pet}), we have
\begin{equation}
\label{normcsphcomp}
\mu_{\normc}^{c_1}(t) \geqslant \mu_{\normc}(t).
\end{equation}

Let us subtract~(\ref{normcF1}) from equality~(\ref{normcF}); then using~(\ref{normcsphcomp}) and the assumption $k_{\normc} \geqslant k_1$ of the theorem, we obtain
\begin{equation}
\label{mainineq2}
\begin{aligned}
0 &\leqslant k_{\normc}(t) - k_1 \\
&\leqslant \frac{d}{dt} \left(\left|\left<N,\partial_{t_O}\right>\right| - \left|\left<N_1,\partial_{t_{O_1}}\right>\right|\right) + \mu_{\normc}^{c_1} (t) \left(\left|\left<N,\partial_{t_O}\right>\right| - \left|\left<N_1,\partial_{t_{O_1}}\right>\right|\right).
\end{aligned}
\end{equation}

If we set $f(t) = \left|\left<N,\partial_{t_O}\right>\right| (t) - \left|\left<N_1,\partial_{t_{O_1}}\right>\right| (t)$, then it follows from~(\ref{mainineq2}) that this function satisfies the following differential inequality 
\begin{equation}
\label{difineq1}
f'(t) + \frac{\sn_{c_1}'(t)}{\sn_{c_1}(t)} f(t) \geqslant 0.
\end{equation}
Since $\sn_{c_1}(t) > 0$ for all positive $t$, inequality~(\ref{difineq1}) is equivalent to
\begin{equation*}
\label{difineq2}
\left( f(t)\cdot\sn_{c_1} (t)\right)' \geqslant 0.
\end{equation*}

Therefore, the function $f\cdot\sn_{c_1}$ is monotonically increasing. Moreover, $f(d)\cdot\sn_{c_1} (d) = 0$. Thus for all $t$ greater then $d$, we have $f(t) \geqslant 0$. Particularly, if $\rho_{O}(P) = \rho_{O_1}(P_1) = l$ ($l > d$), then $f(l) = \left|\left<N,\partial_{t_O}\right>\right| (P) - \left|\left<N_1,\partial_{t_{O_1}}\right>\right| (P_1)\geqslant 0$, as desired.

\begin{remark}
\label{dSrem}
Theorem~\ref{anglecompth1} also holds when $\mathbb M^{m}$ is a de Sitter space $\mathbb S^{m}_1(c)$ of constant positive sectional curvature equal to $c$, $\partial D \subset \mathbb S^{m}_1(c)$ is a connected spacelike hypersurface that is a graph over a standard unit sphere $S^{m-1}$. Such surfaces are called \textit{achronal} (see~\cite{Ger1}).

The assertion above follows from the fact that formula~(\ref{normcF}) can be transferred in the form as it is stated from the Riemannian case to the Lorentzian case almost directly following~\cite{Bor1}. After that one can repeat the calculations from the proof of Theorem~\ref{anglecompth1}.
\end{remark}

\section{Blaschke's Rolling Theorem as a corollary}

In this section we will show that Blaschke's Rolling Theorem is a corollary of Theorems~\ref{anglecompth1} and~\ref{asfcompth}.

We start from the part \textit{A}. Let us introduce in $\mathbb M^{m}(c)$ a polar coordinate system with origin at a point $O \in D$ such that the length of the geodesic segment $OP$ is equal to $\dist(O,\partial D)$. Suppose $(t; \theta^1; \ldots; \theta^{m-1})$ are corresponding coordinates, and assume that the point $P$ has the coordinates $(\dist(O,\partial D); 0; \ldots ;0)$.

Since the domains $D$ and $D_\lambda$ are convex, the hypersurfaces $\partial D$ and $\partial D_\lambda$ that enclose these domains can be given in the introduces coordinate system explicitly by the following equations
\begin{equation}
\label{hypeq}
\partial D: \,t=p(\theta^1,\ldots,\theta^{m-1}),\,\, \partial D_\lambda:\, t = q(\theta^1,\ldots,\theta^{m-1}),
\end{equation}
where $p$ and $q$ are some smooth functions, and $p(0,\ldots,0) = q(0,\ldots,0)$. 

Using~(\ref{hypeq}), we obtain
\begin{equation}
\label{hypangles}
\left|\left<N,\partial_t\right>\right| = \frac{1}{\sqrt{1 + \left|\grad_{\mathbb M}p\right|^2}}, \,\, \left|\left<N_1,\partial_t\right>\right| = \frac{1}{\sqrt{1 + \left|\grad_{\mathbb M}q\right|^2}},
\end{equation}
where $N$ and $N_1$ are inner unit normal fields for, respectively, $\partial D$ and $\partial D_\lambda$; $\partial_t$ is a coordinate vector field tangent to geodesic rays starting from $O$; $\grad_{\mathbb M}$ is a gradient operator defined in $\mathbb M^{m}(c)$. 

If points $Q \in \partial D$, $Q_1 \in \partial D_\lambda$ are such that $\dist(O,Q) = \dist (O,Q_1)$, then by Theorem~\ref{anglecompth1} in a view of~(\ref{hypangles}) at these points the inequality 
\begin{equation*}
\label{mainineq}
\left|\grad_{\mathbb M} p\right| (Q) \leqslant \left|\grad_{\mathbb M} q\right| \left(Q_1\right)
\end{equation*}
holds.

From this point the remaining arguments coincide with those in~\cite[section 4.4]{BorDr2}. And from them it follows that $p \geqslant q$ for all angular parameters $\theta^i$. The last proves the part \textit{A} of Blaschke's Rolling Theorem.

Let us consider the part \textit{B} of the theorem. It is easy to see that for a two-dimensional case ($m=2$) of the part \textit{B} arguments from~\cite{BorDr2} still hold. At the same time, for $m > 2$ they fail to be true. Thus for such a case we need an another approach.

If $\mathbb M^{m}(c)$ is a Euclidean space  $\mathbb E^{m}$, then the part \textit{B} of Blacshke's Rolling Theorem for $m > 2$ follows from the two-dimensional case with a help of projecting. More precisely, if $\pi \subset \mathbb E^{m}$ is an arbitrary two-dimensional plane parallel to a normal vector for $\partial D$ at the point $P$, then an orthogonal projection $Pr_\pi \left(\partial D\right)$ of the hypersurface $\partial D$ on $\pi$ is a curve of curvature not greater than $\lambda$ (see~\cite{Bla} for details).     

If $c \neq 0$, then let us consider a \textit{polar map} of the hypersurface $\partial D$ (see~\cite[Theorem 2.4]{Ger2} and~\cite[Theorem 4.9]{Ger3}). The image of $\partial D$ under this map is a $C^r$-smooth hypersurface $\partial D^\ast$ that lies in a sphere (for $c > 0$), or in a de Sitter space (for $c < 0$). Moreover, normal curvatures $k_{\normc}$ of $\partial D^\ast$ at all points and in every direction satisfy the inequality $k_{\normc} \geqslant 1/\lambda$. Therefore, the hypersurface $\partial D^\ast$ satisfies the part \textit{A} of Blaschke's Rolling Theorem (here we note that, in a view of Remark~\ref{dSrem}, for a de Sitter space all arguments from the proof of the part \textit{A} can be carried out directly). Thus $\partial D^\ast$ lies in a closed convex domain bounded by the hypersurface $\partial D_{1/\lambda}$ of constant normal curvature equal to $1/\lambda$ that touches $\partial D^\ast$ at any given point. Making the polar map of $\partial D^\ast$ and $\partial D_{1/\lambda}$ once more, we will obtain that the complete hypersurface $\partial D_{\lambda} = \left(\partial D_{1/\lambda}\right)^{\ast}$ of constant normal curvature equal to $\lambda$ that touches $\partial D$ at the point $P$ at the same time lies in $D$, as desired. The part \textit{B} is proved.     

\begin{remark}

Blaschke's Rolling Theorem also holds for non-smooth surfaces, namely, when $\partial D$ is a $\lambda$-convex, or $\lambda$-concave hypersurface (for definitions see, for example,~\cite{BorDr2}). This generalized version of Blaschke's Rolling Theorem can be obtained from the smooth version using an approximation result in~\cite[Proposition 6]{PP}.   
\end{remark}

\begin{ack}
This work was partially done while the second author was visiting the Centre de Recerca Matem\`atica as a participant of the Conformal Geometry and Geometric PDE's Program supported by a grant of the Clay Mathematics Institute. He would like to acknowledge both institutions for the given opportunities. 
\end{ack}


\begin{thebibliography}{99}

\bibitem{Bla} W.~Blaschke,
Kreis und Kugel,
~-- de Gruyter, Berlin, 1956.

\bibitem{Kar} {H.~Karcher},
\emph{Umkreise und Inkreise konvexer Kurven in der spharischen und der hyperbolischen Geometrie},
{Math. Ann.}, \textbf{177} (1968), 122-132.

\bibitem{M} {A.D.~Milka}, 
\emph{On a theorem of Schur and Schmidt},
{Ukrain. Geom. Sb.}, \textbf{8} (1970), 95-102. (Russian)

\bibitem{How} {R.~Howard},
\emph{Blaschke's rolling theorem for manifolds with boundary},
{Manuscripta Math.}, \textbf{99} (1999), No. 4, 471-483. 

\bibitem{BZ} Yu.D. {Burago} and V.A. {Zalgaller},
Geometric inequalities (Transl. from Russian by A.B. Sossinsky)
~-- Berlin etc.: Springer-Verlag, 1988.~-- 331 p.

\bibitem{BorDr2} {A.~Borisenko}, {K.~Drach},
\emph{Closeness to spheres of hypersurfaces with normal curvature bounded below},
{Sb. Math.}, \textbf{204} (2013), No. 11, 1565–1583; arXiv:1212.6485.

\bibitem{Bor1} {A.~Borisenko},
\emph{Convex sets in Hadamard manifolds},
Diff. Geom. Appl., \textbf{17} (2002), 111-121.

\bibitem{Pet} \emph{Petersen P.}
\emph{Riemannian geometry},
Graduate texts in mathematics, vol.~171~-- New York: Springer, 1998.

\bibitem{Ger1} C.~Gerhardt,
\emph{Hypersurfaces of prescribed curvature in Lorentzian manifolds},
Indiana U. Math. J., \textbf{49} (2000), No. 3, 1125-1153.

\bibitem{Ger2} C.~Gerhardt,
\emph{Minkowski type problems for convex hypersurfaces in the sphere},
Pure and Applied Mathematics Quarterly, \textbf{3} (2007), No. 2, 417-449.

\bibitem{Ger3} C.~Gerhardt,
\emph{Minkowski type problems for convex hypersurfaces in hyperbolic space},
2006, arXiv:math.DG/0602597, 32 pages.

\bibitem{PP} J.~Parkkonen, F.~Paulin,
\emph{On strictly convex subsets in negatively curved manifolds},
J. Geom. Anal., \textbf{22} (2012), No. 3, 621~-- 632.
\end{thebibliography}
\end{document}